# A self-ruling monotile for aperiodic tiling

Pierre Gradit and Vincent Van Dongen

## Abstract

*Can the entire plane be paved with a single tile that forces aperiodicity? This is known as the ein Stein problem (in German, ein Stein means one tile). This paper presents a monotile that delivers aperiodic tiling by design. It is based on the monotile developed by Taylor and Socolar (whose aperiodicity is forced by means of a non-connected tile that is mainly hexagonal) and motif-based hexagonal tilings that followed this major discovery. Here instead, a single substitution rule makes its shape, and when applying it, forces the tiling to be aperiodic. The proposed monotile, called **HexSeed**, is self-ruling. It consists of 16 identical hexagons, called subtiles, all with edgy borders representing the same binary marking. No motif is needed on the subtiles to make it work. Additional motifs can be added to the monotile to provide some insights. The proof of aperiodicity is presented with the use of such motifs.*

## Introduction

In the seventies, Penrose demonstrated that two aperiodic tiles can pave the plane without holes. The two tiles are said to be aperiodic because they can only generate an aperiodic tiling. Soon after, people wondered whether an aperiodic monotile could exist. That is, could a single tile be aperiodic, not allowing any periodic tiling, and still pave the entire plane? This problem is known as the *'ein-Stein'* problem. ('ein Stein' is German and means 'one tile'.)

The existence of a single connected tile to pave the entire plane only in an aperiodic manner with no other artefact than its shape is still unknown. However, some major discoveries were made. A decade ago, Taylor and Socolar (Socolar, 2011) proposed a first solution to this problem. The tile is mainly hexagonal but with non-connected borders to enforce the aperiodicity. It is shown here below as well as the tiling it generates (referred hereafter as T-S tiling).

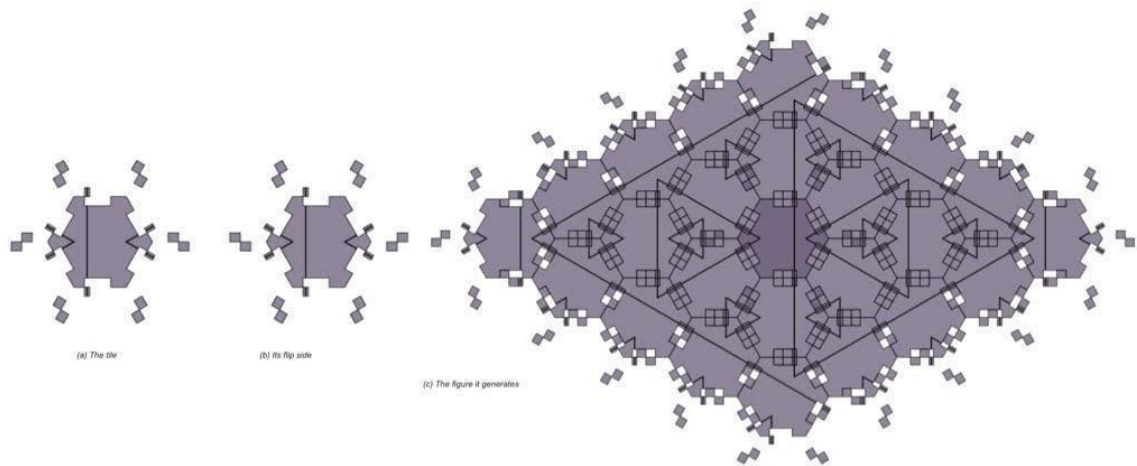

*Figure 1: Non-connected tile proposed as a solution to the ein-Stein problem and the Sierpinski triangle it generates when using it. Ref: https://en.wikipedia.org/wiki/Einstein_problem*

Note that the motif on the tile (the black lines) is not necessary as the shape of the tile alone creates the aperiodicity. There is only one way of covering the plane with it and it requires both sides of the tile (the other side is obtained by flipping the tile). The black lines on the tile provide some insight on the tiling being generated. On the paved area, a figure gets created that is a fractal similar to the well-known triangle of Sierpinski (Sierpinski triangle). See figure here above on the right.

In T-S tiling, the non-connected flags allow the constraints needed for aperiodicity to be non strictly local. In the same publication (Socolar, 2011), a connected version of their monotile is also presented. This connected version is of hexagonal shape with a motif for forcing the aperiodicity. Note that the motif imposes constraints on non-adjacent tiles.

In (Penrose, 1997), another hexagonal-based tiling, called $(1 + \varepsilon + \varepsilon^2)$-tiling, is presented with thin edge tiles and small corner tiles. Both T-S tiling and Penrose tiling are compared in (Baake, Gahler, & Grimm, Hexagonal Inflation Tilings and Planar Monotiles, 2012) and (Lee & Moody, 2017). It is noted that both tilings are based on a hierarchical system of nested equilateral triangles, i.e. the Sierpinski triangle.

An alternative motif on hexagonal tiles was recently proposed in (Mampusti & Whittaker, 2020). The solution makes use of a dendrite motif. As explained in their paper, the dendrite forms nested triangles of Sierpinski as well. The innovation here is that the monotile does not require to be flipped. The enforcement of aperiodicity is controlled by the dendrite property. This dendrite method will be the basis of our study on aperiodic tiling.

In this paper, we present yet another version of a hexagonal based tiling called HexSeed. Its aperiodicity is enforced by making our monotile self-ruling as explained in the next section.

## A self-ruling monotile called HexSeed

A self-ruling monotile is a single tile self-containing scale and substitution information. Our monotile called HexSeed is shown here below.

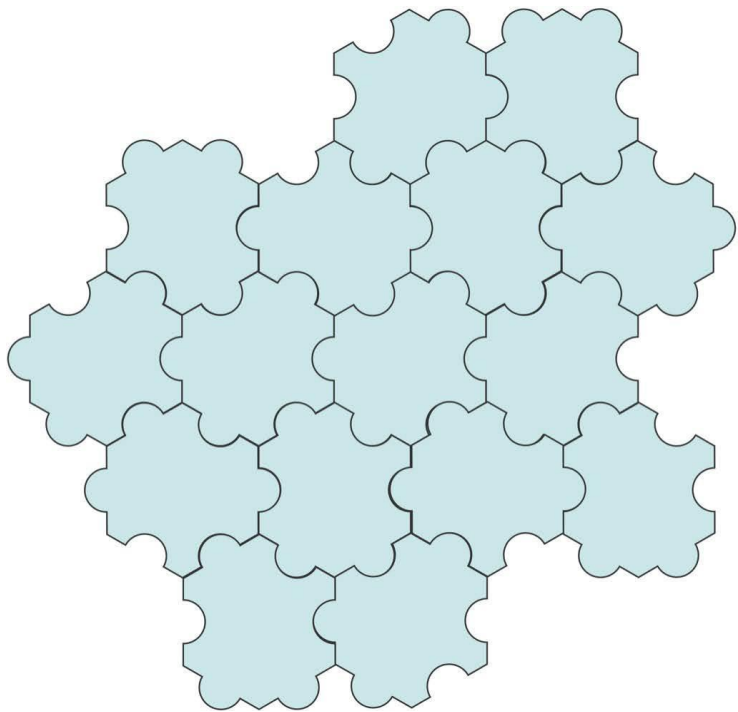

*Figure 2: HexSeed is a self-ruling monotile as it self-contains scale and substitution information.*

The scaling and substitution information contained in HexSeed are to be used to develop the aperiodic tiling. In other words, HexSeed is to be built iteratively out of the scaling/substitution rule shown here below.

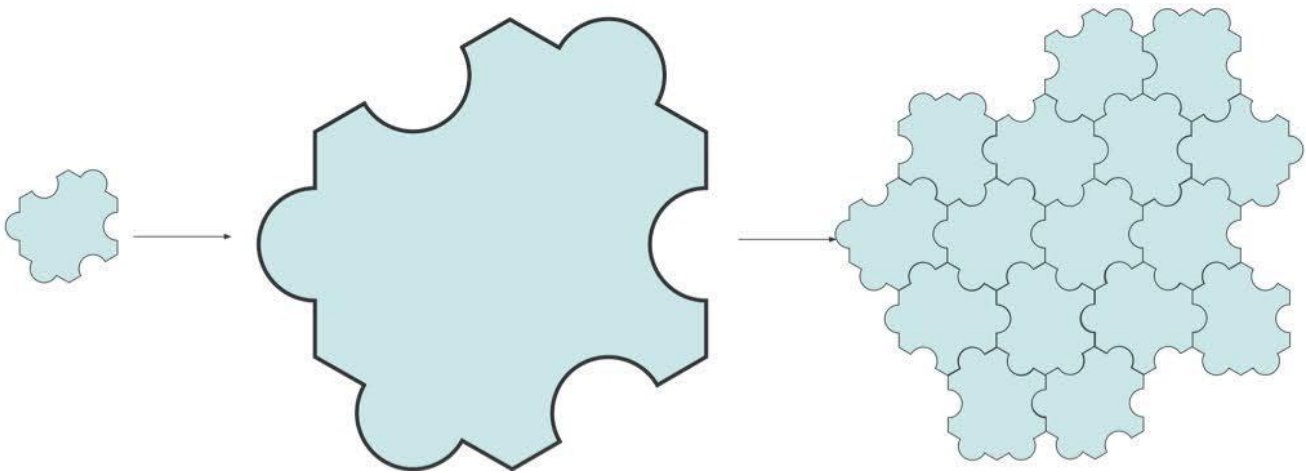

*Figure 3: HexSeed Scale/Substitution rule as depicted by its right-hand side alone.*

The monotile HexSeed consists of 16 identical subtiles and the substitution rule defines the exact orientation of each of them. Each subtile is a motif-free hexagon with borders marked by their shape.

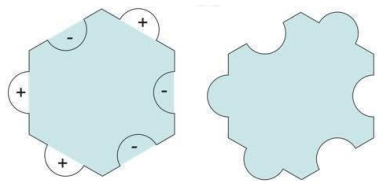

*Figure 4: HexSeed subtiles with bumps (+) and holes (-)*

These borders are shaped in such a way that, once enlarged, both subtile and monotile closely match, as shown here below.

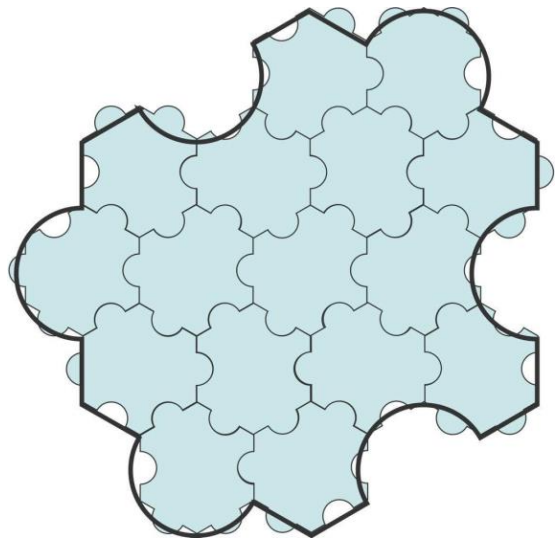

*Figure 5: The close match between the enlarged subtile and the monotile.*

Note that **HexSeed** can come in diffe5ent versions. Its subtile is an hexagon with edgy borders that need to look, once enlarged to the monotile as close as possible in order to ease its use.

**HexSeed** is closely related to hexagon-based aperiodic tilings presented in the introduction. In order to demonstrate this, we need to add motifs to the tiles. In the remaining of this paper, we will consider two ways of applying the scale/substitution rule:

- **the finite way**, making use of a finite number of application of the rule
- **the infinite way**, by considering a plan partition which is a fix point of the rule

# A finite exploration with motifs

In order to explore the properties of the tile, we first propose to assign a single bit code, for either black or white, to each subtile. This already provides many possibilities and offers the best contrasting result. Using this exploratory technique, two main patterns are recurrent: Sierpinski triangle and the dendrite. First pattern, well identified in literature, as mentioned in the introduction, is the Sierpinski triangle. See on the figure here below how it can be generated with a single line of black subtiles on the monotile.

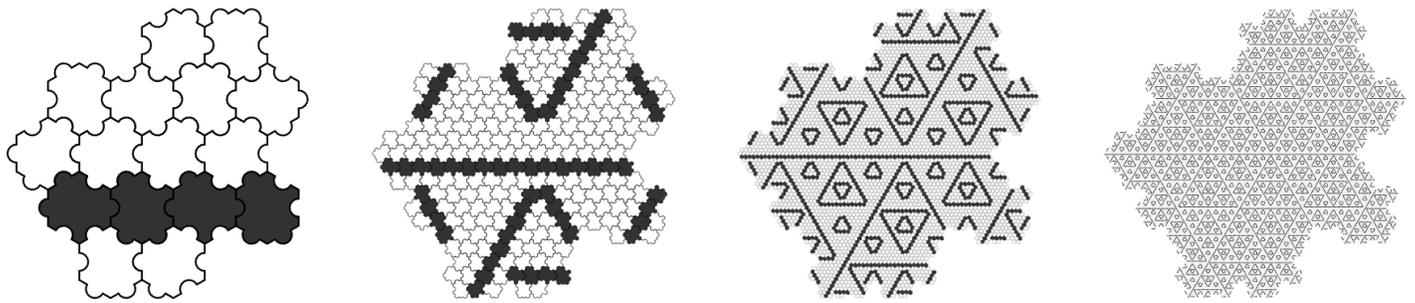

*Figure 6: The Sierpinski triangle can be generated by assigning a single line of subtiles in black on the monotile. Here the substitution rule is applied three times (without scaling, hence all images are of the same size).*

The other motif of particular interest is the dendrite. It can be created as shown here below.

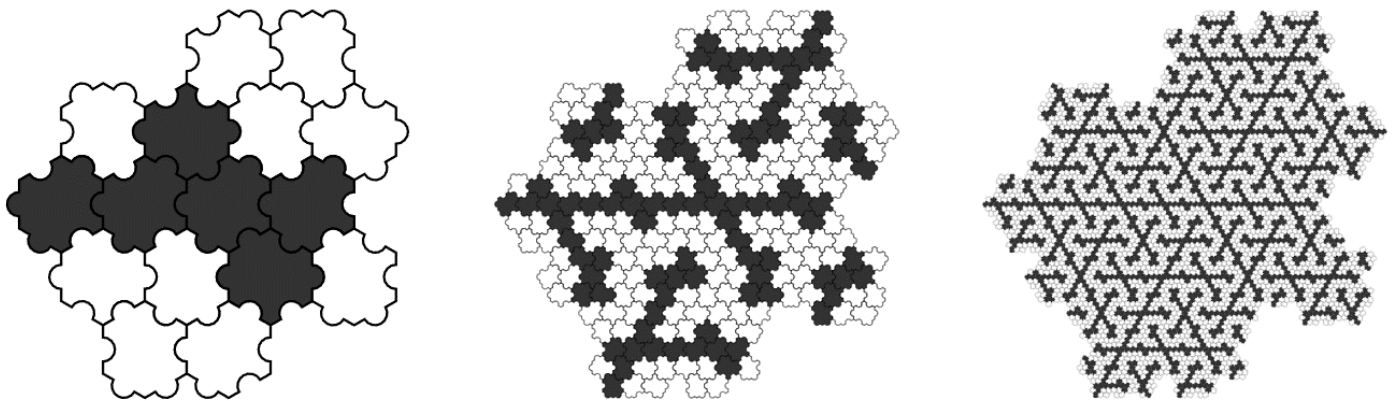

*Figure 7: The dendritic motif in black and white on the monotile and two applications of its substitution rule.*

Let us now make use of a similar motif but with finer details on the subtile in order to better analyse the aperiodic property of **HexSeed**. With the representation provided on the figure here below, dendrites convey the idea of infinite paths originating from sources represented by circles.

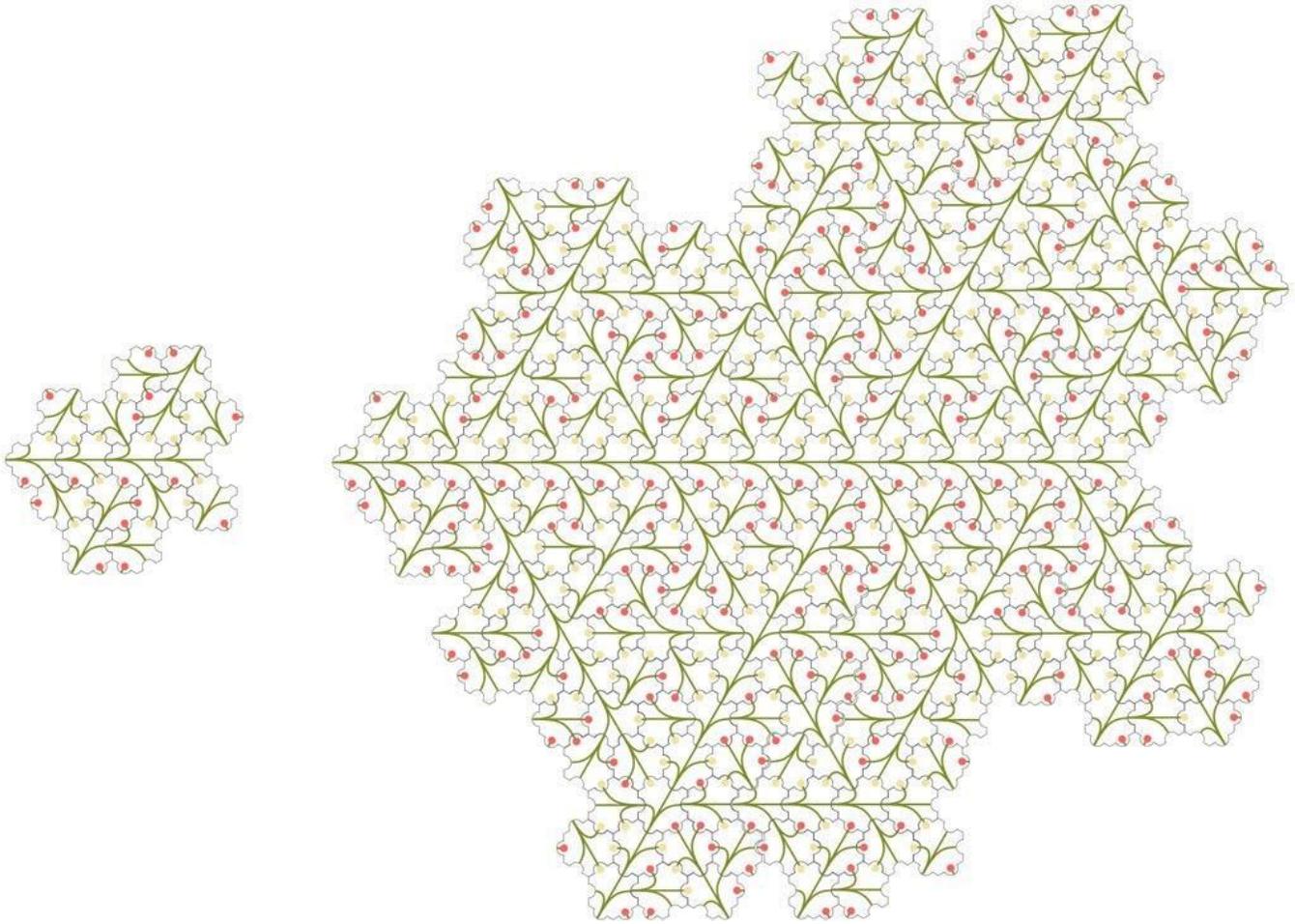

*Figure 8: Monotile with motif of dendrite with sources represented by circles.*

In 2020, it was proven in (Mampusti & Whittaker, 2020) that if dendrites have no cycle, the tiling is aperiodic. In the next section, we will indeed prove that the dendrites generated by **HexSeed** contain no cycle.

It can be observed that the subtiles in HexSeed take six orientations. For clarity, we have assigned a different color for each orientation. The tiling can be created in two steps. First, an infinite star that divides the plane in 6 is created as shown here below on the left. This property, given at the left of Fig 9, is necessary to initiate our final proof in the next part.

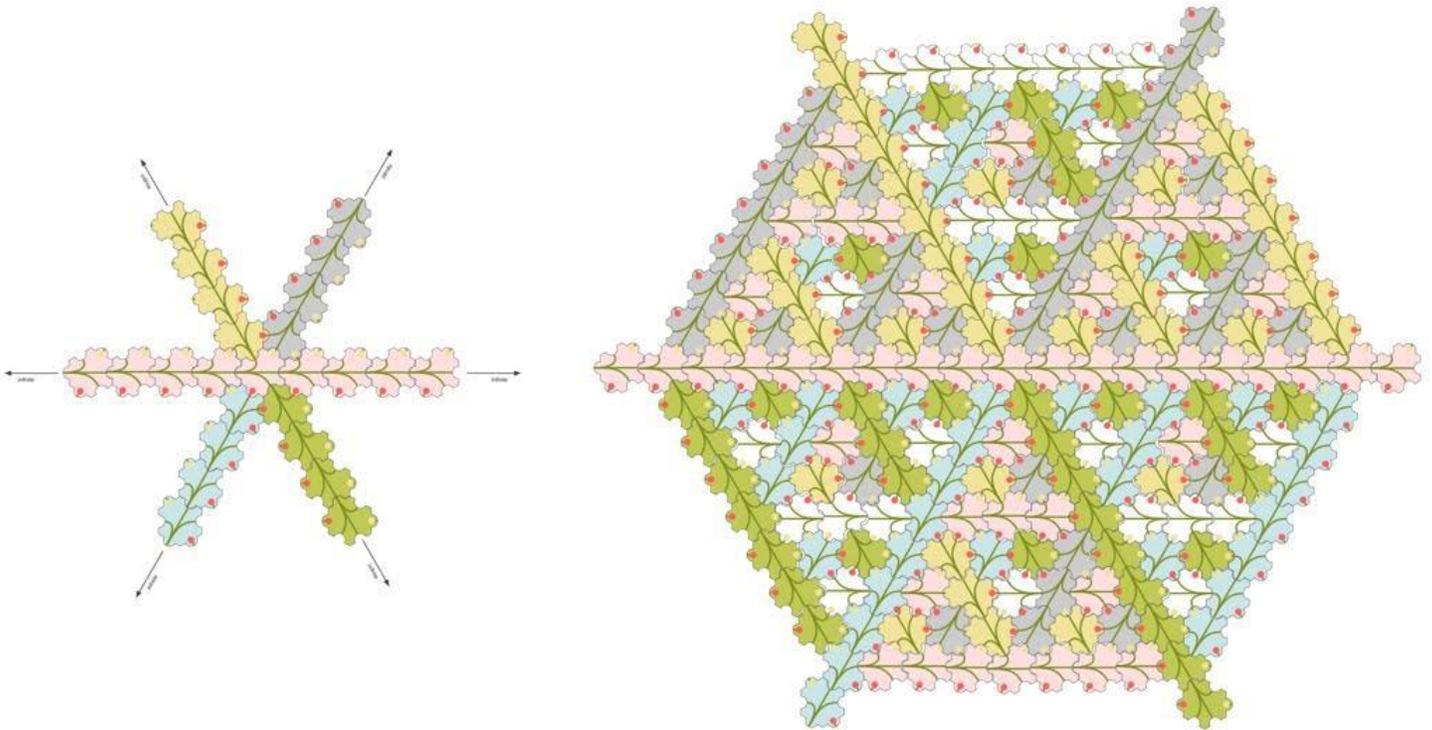

*Figure 9: Tiling identical to HexSeed can be created directly from the subtiles in two steps: First (see on the left) an infinite star is created that divides the plane in 6 and then (see on the right) the tiling is completed starting from the center of the star.*

The plane then can be tiled from the center of the star. The subtiles can simply be placed as imposed by the star through a process we will call "must-hex" in the latter. In order to prove the aperiodicity through that path, we would need to state that the single colored lines have a pattern as we can grasp it visually on the right-hand side of Figure 9. But this way to make the proof, triangle based, was not so straight forward and we choose instead the dendrite approach presented here below.

## Proving the aperiodicity of HexSeed

In order to prove with sufficient safety that our monotile scale/substitution rule does not create dendrite cycles, we need to consider infinite scale/substitution rule application by calling P the infinite plane partition of **HexSeed** which is the fix-point of the scale/substitution rule application.

Our proof of aperiodicity relies on the capacity to colour all sources: white if they go to the -180° dendrite and black if they go to either the 60° dendrite or the -120° depending on their position with respect to the horizontal infinite dendrite. If that colouring is decidable then any source will go into an infinite dendrite, which implies it has no cycle, and therefore the tiling is aperiodic.

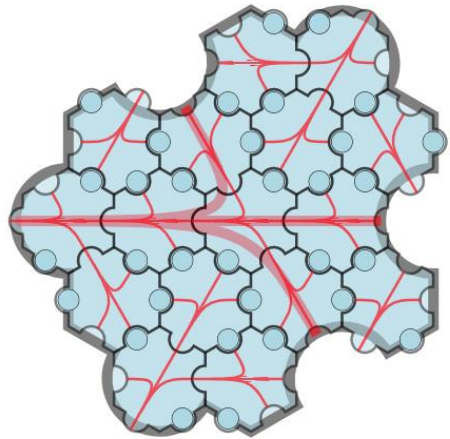

*Figure 10: The tile with local dendrites.*

**The first stage** of the proof is to complete the rule single application with all "must-hex", i.e. subtiles that must be placed that way because of the occurrence of a concave sequence of similar shapes: two bumps (++) or two holes (--). This completion stage leads to the following figure after three iterations (green, yellow, orange), where all concavities are alternate bump/hole sequence:

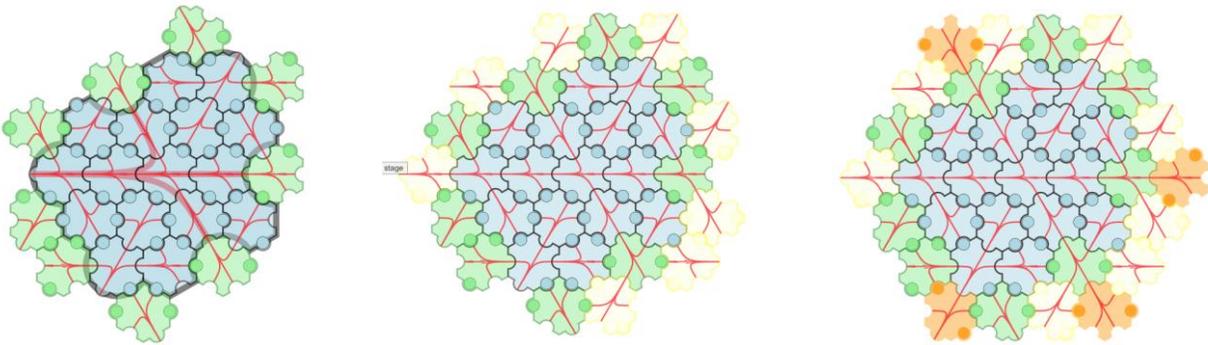

*Figure 11: completing the HexSeed tile with "must-hexes".*

**The second stage** of the proof is a diffusion computing of the knowledge about black/white sources. As seen in Figure 10, the center subtile also has the property of being connected to all dendrites as it is situated on the -180° horizontal dendrite and contains the sources of both 60° and -120° dendrites. The initiation of the computing being that the origin subtile is on the infinite -180° dendrite and both its sources are initiation of respective black dendrites.

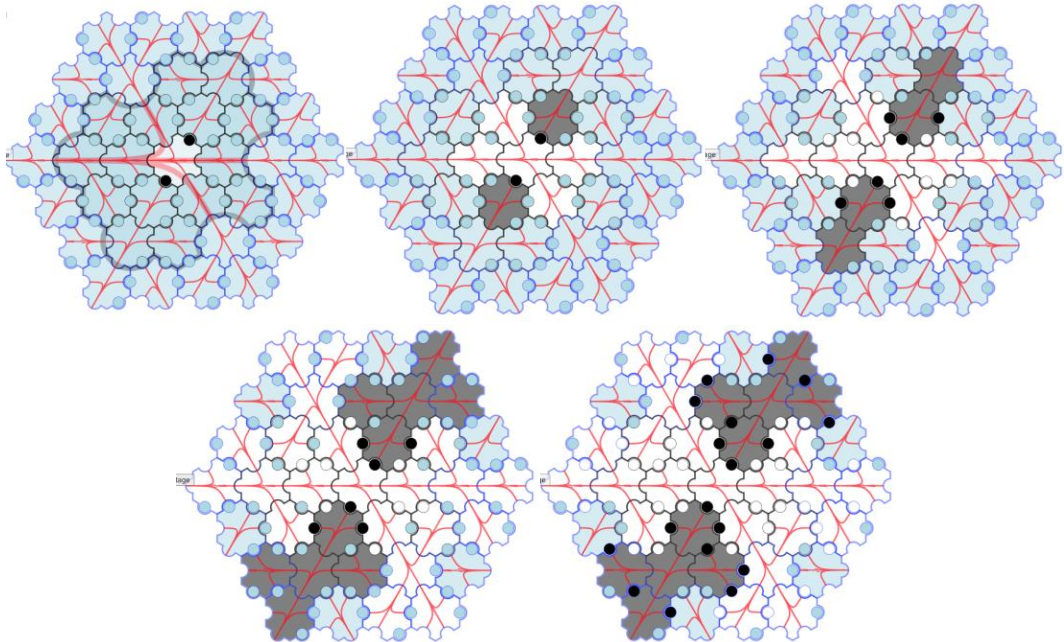

*Figure 12: Deciding colour by diffusion leaving six undecided subtiles.*

**The third stage** of the proof is to decide the colour of the sources of the six blue subtiles on the border by applying the scale/substitution rule.

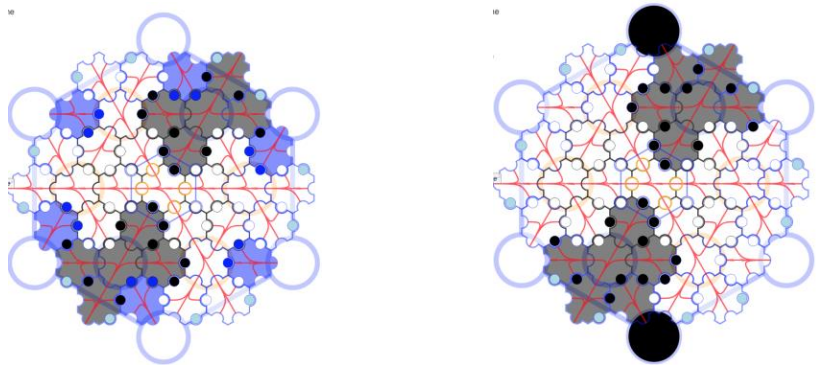

*Figure 13: Use of scale/substitution information to conclude*

**[Aperiodicity proof]** "*a source is either 180°, 60° or -120°*" is a decidable predicate.

The last stage of the proof, illustrated in Fig. 13, relies on that infinite rule application that allow to state on the perimeter sources of the tile. As soon as the color of any source of the self-ruling tile is decided, any source of any subtile of the plane can be decided with the same algorithm.

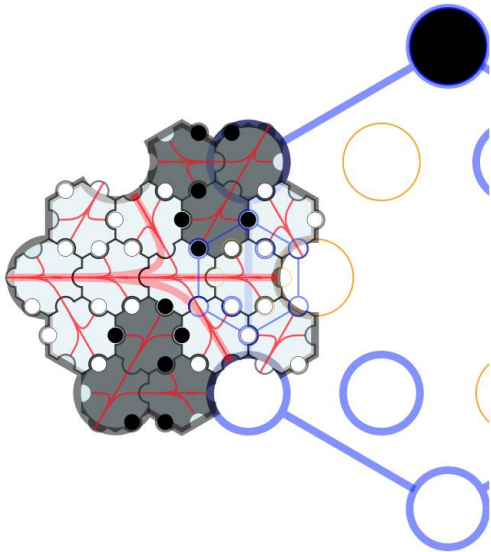

*Figure 14: Walking in the infinite plane*

You can walk wherever you want on the plane which is the fix point and compute the color of all sources on your path. For example, walking at 0° as in figure 15, you consider the valuation to apply for the indicated expanded hexagon. Its substitution will be composed of white sources except three perimeter sources on the upper side, black. These three sources are the two at the top and the one at upper-left edge.

This process can be written by identifying each subtile with a hexadecimal, i.e. any subtile tiling the plane can be associated with a hexadecimal word. We can then build an automaton conveying orientation and the twelve defining values to compute the source colour for any source of the plane. There is no major difficulty in this computing once you get the way to decide the color of the perimeter sources by the completion and fix point hypothesis.

# Conclusion

The study of the "ein-stein" problem as opened in 2010 by Taylor and Socolar seems to be always associated with an hexagonal base and a rule conveying a non-strictly adjacent rule to insure aperiodicity. Our work embeds this constraint in a single tile, called **HexSeed**, that conveys sufficient information to tile the whole plane aperiodically and with no tile flipping.

By colouring its subtiles, we made apparent two motifs that are useful for proving its aperiodicity: the Sierpinski Triangle and the dendrite motif. In this paper we choose the dendrite motif and prove aperiodicity by deciding the final and infinite destination of any source (-180°, 60° or -120°). As these dendrites are infinite, the absence of cycle and thus aperiodicity is proven.

We showed in this paper the importance of selecting the right motif for the study of aperiodic tilings. We feel this exploration is far from being completed and therefore plan to pursue it in the future. Other ways to convey the non-local constraint may be found, even lighter than these we know as granted, or the one we propose in this paper.

Behind the potential applications that these tilings may bring, let us not forget the aesthetic value of them as it is how it all started. The following figure is a particular example of **HexSeed** with a motif that we particularly found pleasant to look at. This invaluable aspect remains an important part of our interests in tilings.

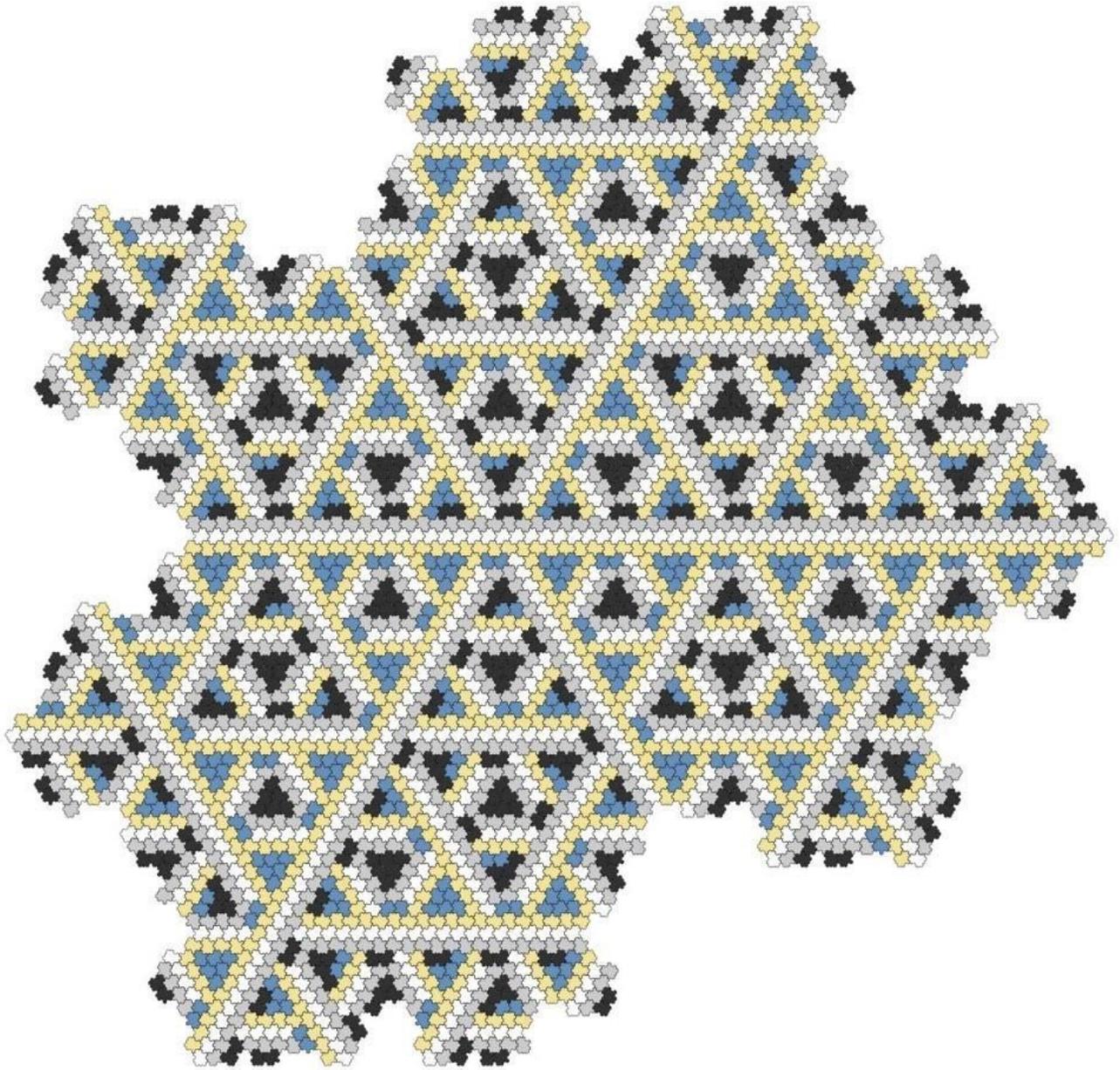

*Figure 15: An example of **HexSeed** with a motif that we particularly enjoyed looking at from a certain distance.*